\documentclass{gtart}
\def\ifplaintex{\expandafter\ifx\csname documentclass\endcsname\relax}

\def\gtp{{\mathsurround=0pt\it $\cal G\mskip-2mu$eometry \&\ 
$\cal T\!\!$opology $\cal P\!$ublications}}  % GT publications

\def\recd{{\small Received:\qua\receiveddate\ifx\reviseddate\relax
\else\qquad Revised:\qua\reviseddate\fi\par}} 

\def\lognumber#1{\def\thelognumber{#1}}
\def\volumenumber#1{\def\thevolumenumber{#1}}
\def\volumeyear#1{\def\thevolumeyear{#1}}
\def\papernumber#1{\def\thepapernumber{#1}}
\def\pagenumbers#1#2{\def\startpage{#1}\def\finishpage{#2}}
\def\published#1{\def\publishdate{#1}}

\def\received#1{\def\receiveddate{#1}}
\def\revised#1{\def\reviseddate{#1}}
\def\accepted#1{\def\accepteddate{#1}}
\def\asciititle#1{\def\theasciititle{#1}}
\def\covertitle#1{\def\thecovertitle{#1}}

\def\asciiemail#1{\def\theasciiemail{#1}}

\long\def\asciiabstract#1{\long\def\theasciiabstract{#1}}
\def\asciikeywords#1{\def\theasciikeywords{#1}}

\let\\\par\let\thelognumber\relax\let\thevolumenumber\relax
\let\thepapernumber\relax\let\thevolumeyear\relax\let\startpage\relax
\let\finishpage\relax\let\publishdate\relax\let\receiveddate\relax
\let\reviseddate\relax\let\accepteddate\relax\let\theasciititle\relax
\let\thecovertitle\relax\let\theasciiauthors\relax
\let\theasciiabstract\relax\let\theasciikeywords\relax

\let\theasciiemail\relax

\ifplaintex
\font\logobig=cmssbx10 scaled 3836
\font\logomed=cmssbx10 scaled 2557
\else
\font\logobig=cmssbx10 scaled 4200
\font\logomed=cmssbx10 scaled 2800
\fi

\long\def\makeagttitle{   %%% start of definition of \makeagttitle
\count0=\startpage
\agt\hfill      %   Journal title (top left) 
\hbox to 45truept{\vbox to 0pt{\vglue -13truept{\logomed A\kern -.37em{\logobig 
T}\kern -.38em G}\vss}\hss}
\break
{\small Volume \thevolumenumber\ (\thevolumeyear)
\startpage--\finishpage\nl
Published: \publishdate}

\vglue .25truein

{\parskip=0pt\leftskip 0pt plus
1fil\def\\{\par\smallskip}{\Large\bf\thetitle}\par\medskip} \vglue
0.05truein

{\parskip=0pt\leftskip 0pt plus 1fil\def\\{\par}{\sc\theauthors}
\par\medskip}%
 
\vglue 0.03truein 

{\small\leftskip 25truept\rightskip 25truept{\bf Abstract}\stdspace\theabstract

{\bf AMS Classification}\stdspace\theprimaryclass
\ifx\thesecondaryclass\relax\else; \thesecondaryclass\fi\par
{\bf Keywords}\stdspace \thekeywords\par}\vglue 7truept

}   %%%% end of definition of \makeagttitle

\ifplaintex
\font\phead=cmsl9 scaled 950
\font\pnum=cmbx10 scaled 913
\font\pfoot=cmsl9 scaled 950
\headline{\vbox to 0pt{\vskip -4.5mm\line{\small\phead\ifnum
\count0=\startpage ISSN 1472-2739 (on-line) 1472-2747 (printed)
\hfill {\pnum\folio}\else\ifodd\count0\def\\{ }% 
\ifx\theshorttitle\relax\thetitle\else\theshorttitle\fi\hfill{\pnum\folio}
\else\def\\{ and }{\pnum\folio}\hfill\ifx\theshortauthors\relax\theauthors
\else\theshortauthors\fi\fi\fi}\vss}}
\footline{\vbox to 0pt{\vglue 0mm\line{\small\pfoot\ifnum\count0=\startpage
\copyright\ \gtp\hfill\else
\agt, Volume \thevolumenumber\ (\thevolumeyear)\hfill\fi}\vss}}
\else
\font=cmsl9 scaled 1050
\font\lnum=cmbx10 
\font=cmsl9 scaled 1050
\makeatletter
\def\@oddhead{{\small\lhead\ifnum\count0=\startpage ISSN 1472-2739 
(on-line) 1472-2747 (printed)\hfill {\lnum\number\count0}\else\ifodd\count0
\def\\{ }\ifx\theshorttitle\relax \thetitle \else\theshorttitle\fi\hfill
{\lnum\number\count0}\else\def\\{ and }{\lnum\number\count0}
\hfill\ifx\theshortauthors\relax 
\theauthors\else\theshortauthors\fi\fi\fi}}\def\@evenhead{\@oddhead}
\def\@oddfoot{\small\lfoot\ifnum\count0=\startpage\copyright\ \gtp\hfill\else
\agt, Volume \thevolumenumber\ (\thevolumeyear)\hfill\fi}
\def\@evenfoot{\@oddfoot}
\makeatother
\fi
\let\maketitlepage\makeagttitle
\let\makeshorttitle\maketitlepage
\let\maketitle\maketitlepage

\newwrite\gtoutfile
\long\gdef\makeheadfile{  %%% start of definition of \makeheadfile
{\def\\{, }\def\s{ }
\immediate\openout\gtoutfile head.xxx
\immediate\write\gtoutfile{To: math@arxiv.org}
\immediate\write\gtoutfile{Subject: put OR rep NNNNN:ppppp}
\immediate\write\gtoutfile{--text follows this line--}
\immediate\write\gtoutfile{Proxy-for: \ifx\theasciiauthors\relax
\theauthors\else\theasciiauthors\fi\s<\ifx\theasciiemail\relax\theemail\else\theasciiemail\fi>}
\immediate\write\gtoutfile{\noexpand\\}
\immediate\write\gtoutfile{Authors: \ifx\theasciiauthors\relax
\theauthors\else\theasciiauthors\fi}
{\def\\{ }\immediate\write\gtoutfile{Title: \ifx\theasciititle\relax
\thetitle\else\theasciititle\fi}}
\immediate\write\gtoutfile{Subj-class: GT or SG, GR etc}
\immediate\write\gtoutfile{MSC-class: \theprimaryclass\ifx\thesecondaryclass\relax\else, \thesecondaryclass\fi}
\immediate\write\gtoutfile{Journal-ref: Algebr. Geom. Topol. \thevolumenumber\s
(\thevolumeyear) \startpage-\finishpage}
\immediate\write\gtoutfile{Comments: Published by Algebraic and
Geometric Topology at}
\immediate\write\gtoutfile{\s\s\s  http://www.maths.warwick.ac.uk/agt/AGTVol\thevolumenumber/agt-\thevolumenumber-\thepapernumber.abs.html}
\immediate\write\gtoutfile{\noexpand\\}
\immediate\write\gtoutfile{}
\ifx\theasciiabstract\relax
\immediate\write\gtoutfile{\theabstract}\else
\immediate\write\gtoutfile{\theasciiabstract}\fi
\immediate\write\gtoutfile{}
\immediate\write\gtoutfile{\noexpand\\}
\immediate\write\gtoutfile{}
\immediate\closeout\gtoutfile}}  %%% end of definition of \makeheadfile

\def\maketitlepage{\makeagttitle\makeheadfile}
\let\makeshorttitle\maketitlepage
\let\maketitle\maketitlepage

\def\ifplaintex{\expandafter\ifx\csname documentclass\endcsname\relax}

\def\gtp{{\mathsurround=0pt\it $\cal G\mskip-2mu$eometry \&\ 
$\cal T\!\!$opology $\cal P\!$ublications}}  % GT publications

\def\recd{{\small Received:\qua\receiveddate\ifx\reviseddate\relax
\else\qquad Revised:\qua\reviseddate\fi\par}} 

\def\lognumber#1{\def\thelognumber{#1}}
\def\volumenumber#1{\def\thevolumenumber{#1}}
\def\volumeyear#1{\def\thevolumeyear{#1}}
\def\papernumber#1{\def\thepapernumber{#1}}
\def\pagenumbers#1#2{\def\startpage{#1}\def\finishpage{#2}}
\def\published#1{\def\publishdate{#1}}

\def\received#1{\def\receiveddate{#1}}
\def\revised#1{\def\reviseddate{#1}}
\def\accepted#1{\def\accepteddate{#1}}
\def\asciititle#1{\def\theasciititle{#1}}
\def\covertitle#1{\def\thecovertitle{#1}}

\def\asciiemail#1{\def\theasciiemail{#1}}

\long\def\asciiabstract#1{\long\def\theasciiabstract{#1}}
\def\asciikeywords#1{\def\theasciikeywords{#1}}

\let\\\par\let\thelognumber\relax\let\thevolumenumber\relax
\let\thepapernumber\relax\let\thevolumeyear\relax\let\startpage\relax
\let\finishpage\relax\let\publishdate\relax\let\receiveddate\relax
\let\reviseddate\relax\let\accepteddate\relax\let\theasciititle\relax
\let\thecovertitle\relax\let\theasciiauthors\relax
\let\theasciiabstract\relax\let\theasciikeywords\relax

\let\theasciiemail\relax

\ifplaintex
\font\logobig=cmssbx10 scaled 3836
\font\logomed=cmssbx10 scaled 2557
\else
\font\logobig=cmssbx10 scaled 4200
\font\logomed=cmssbx10 scaled 2800
\fi

\long\def\makeagttitle{   %%% start of definition of \makeagttitle
\count0=\startpage
\agt\hfill      %   Journal title (top left) 
\hbox to 45truept{\vbox to 0pt{\vglue -13truept{\logomed A\kern -.37em{\logobig 
T}\kern -.38em G}\vss}\hss}
\break
{\small Volume \thevolumenumber\ (\thevolumeyear)
\startpage--\finishpage\nl
Published: \publishdate}

\vglue .25truein

{\parskip=0pt\leftskip 0pt plus
1fil\def\\{\par\smallskip}{\Large\bf\thetitle}\par\medskip} \vglue
0.05truein

{\parskip=0pt\leftskip 0pt plus 1fil\def\\{\par}{\sc\theauthors}
\par\medskip}%
 
\vglue 0.03truein 

{\small\leftskip 25truept\rightskip 25truept{\bf Abstract}\stdspace\theabstract

{\bf AMS Classification}\stdspace\theprimaryclass
\ifx\thesecondaryclass\relax\else; \thesecondaryclass\fi\par
{\bf Keywords}\stdspace \thekeywords\par}\vglue 7truept

}   %%%% end of definition of \makeagttitle

\ifplaintex
\font\phead=cmsl9 scaled 950
\font\pnum=cmbx10 scaled 913
\font\pfoot=cmsl9 scaled 950
\headline{\vbox to 0pt{\vskip -4.5mm\line{\small\phead\ifnum
\count0=\startpage ISSN 1472-2739 (on-line) 1472-2747 (printed)
\hfill {\pnum\folio}\else\ifodd\count0\def\\{ }% 
\ifx\theshorttitle\relax\thetitle\else\theshorttitle\fi\hfill{\pnum\folio}
\else\def\\{ and }{\pnum\folio}\hfill\ifx\theshortauthors\relax\theauthors
\else\theshortauthors\fi\fi\fi}\vss}}
\footline{\vbox to 0pt{\vglue 0mm\line{\small\pfoot\ifnum\count0=\startpage
\copyright\ \gtp\hfill\else
\agt, Volume \thevolumenumber\ (\thevolumeyear)\hfill\fi}\vss}}
\else
\font=cmsl9 scaled 1050
\font\lnum=cmbx10 
\font=cmsl9 scaled 1050
\makeatletter
\def\@oddhead{{\small\lhead\ifnum\count0=\startpage ISSN 1472-2739 
(on-line) 1472-2747 (printed)\hfill {\lnum\number\count0}\else\ifodd\count0
\def\\{ }\ifx\theshorttitle\relax \thetitle \else\theshorttitle\fi\hfill
{\lnum\number\count0}\else\def\\{ and }{\lnum\number\count0}
\hfill\ifx\theshortauthors\relax 
\theauthors\else\theshortauthors\fi\fi\fi}}\def\@evenhead{\@oddhead}
\def\@oddfoot{\small\lfoot\ifnum\count0=\startpage\copyright\ \gtp\hfill\else
\agt, Volume \thevolumenumber\ (\thevolumeyear)\hfill\fi}
\def\@evenfoot{\@oddfoot}
\makeatother
\fi
\let\maketitlepage\makeagttitle
\let\makeshorttitle\maketitlepage
\let\maketitle\maketitlepage

\newwrite\gtoutfile
\long\gdef\makeheadfile{  %%% start of definition of \makeheadfile
{\def\\{, }\def\s{ }
\immediate\openout\gtoutfile head.xxx
\immediate\write\gtoutfile{To: math@arxiv.org}
\immediate\write\gtoutfile{Subject: put OR rep NNNNN:ppppp}
\immediate\write\gtoutfile{--text follows this line--}
\immediate\write\gtoutfile{Proxy-for: \ifx\theasciiauthors\relax
\theauthors\else\theasciiauthors\fi\s<\ifx\theasciiemail\relax\theemail\else\theasciiemail\fi>}
\immediate\write\gtoutfile{\noexpand\\}
\immediate\write\gtoutfile{Authors: \ifx\theasciiauthors\relax
\theauthors\else\theasciiauthors\fi}
{\def\\{ }\immediate\write\gtoutfile{Title: \ifx\theasciititle\relax
\thetitle\else\theasciititle\fi}}
\immediate\write\gtoutfile{Subj-class: GT or SG, GR etc}
\immediate\write\gtoutfile{MSC-class: \theprimaryclass\ifx\thesecondaryclass\relax\else, \thesecondaryclass\fi}
\immediate\write\gtoutfile{Journal-ref: Algebr. Geom. Topol. \thevolumenumber\s
(\thevolumeyear) \startpage-\finishpage}
\immediate\write\gtoutfile{Comments: Published by Algebraic and
Geometric Topology at}
\immediate\write\gtoutfile{\s\s\s  http://www.maths.warwick.ac.uk/agt/AGTVol\thevolumenumber/agt-\thevolumenumber-\thepapernumber.abs.html}
\immediate\write\gtoutfile{\noexpand\\}
\immediate\write\gtoutfile{}
\ifx\theasciiabstract\relax
\immediate\write\gtoutfile{\theabstract}\else
\immediate\write\gtoutfile{\theasciiabstract}\fi
\immediate\write\gtoutfile{}
\immediate\write\gtoutfile{\noexpand\\}
\immediate\write\gtoutfile{}
\immediate\closeout\gtoutfile}}  %%% end of definition of \makeheadfile

\def\maketitlepage{\makeagttitle\makeheadfile}
\let\makeshorttitle\maketitlepage
\let\maketitle\maketitlepage

\lognumber{9}
\volumenumber{1}
\volumeyear{2001}
\papernumber{9}
\pagenumbers{173}{199}
\received{19 February 2001}
\revised{4 April 2001}
\accepted{6 April 2001}
\published{7 April 2001}

\let\cal\mathcal

\usepackage{amssymb,amscd,amsxtra,eucal}
\usepackage[all]{xy}

\newtheorem{thm}{Theorem}[section]
\newtheorem{lem}[thm]{Lemma}
\newtheorem{prop}[thm]{Proposition}
\newtheorem{cor}[thm]{Corollary}

\theoremstyle{definition}
\newtheorem{rem}[thm]{Remark}

\newtheorem{conj}[thm]{Conjecture}
\newtheorem{cond}[thm]{Conditions}

\numberwithin{equation}{section}

\def\ds{\displaystyle}
\def\:{\colon}
\def\.{\cdot}

\def\<{\left\langle}
\def\>{\right\rangle}
\def\({\left(}
\def\){\right)}
\def\ph#1{\phantom{#1}}
\def\epsilon{\varepsilon}
\def\phi{\varphi}

\def\leq{\leqslant}
\def\geq{\geqslant}
\def\la{\leftarrow}
\def\ra{\rightarrow}
\def\lla{\longleftarrow}
\def\lra{\longrightarrow}
\def\Lra{\Longrightarrow}
\def\bar#1{\overline{#1}}
\def\hat#1{\widehat{#1}}
\def\tilde#1{\widetilde{#1}}
\def\iso{\cong}
\DeclareMathOperator{\DERIV}{d}\renewcommand{\d}{\DERIV}

\DeclareMathOperator{\cofibre}{cofibre}

\DeclareMathOperator{\fibre}{fibre}

\DeclareMathOperator{\coker}{coker}

\def\E{\mathrm{E}}
\def\F{\mathbb{F}}

\def\Smash{\wedge}
\def\SSmash#1{{\ds\mathop{\Smash}_{#1}}}
\def\Square#1{\mathop{\square}_{#1}}

\def\oTimes#1{{\ds\mathop{\otimes}_{#1}}}

\def\ideal{\triangleleft}
\DeclareMathOperator{\Coext}{Coext}
\DeclareMathOperator{\Cotor}{Cotor}

\DeclareMathOperator{\Ext}{Ext}
\DeclareMathOperator{\Func}{F}
\DeclareMathOperator{\Hom}{Hom}

\DeclareMathOperator{\Tor}{Tor}
\def\BP{{BP}}
\def\En{E(n)}
\def\Kn{K(n)}

\def\MSp{MSp}
\def\MU{{MU}}
\def\MUn{\MU\<n\>}
\DeclareMathOperator*{\colim}{colim}
\DeclareMathOperator*{\invlim}{lim}
\DeclareMathOperator*{\holim}{holim}
\DeclareMathOperator*{\hocolim}{hocolim}

\def\Spectra{{\mathcal S}}

\def\Func{\mathcal F}

\DeclareMathOperator{\LOC}{L}\renewcommand{\L}{\LOC}
\def\del{\partial}
\DeclareMathOperator{\bideg}{bideg}
\begin{document}
\title[On the Adams Spectral Sequence for $R$-modules]{On the Adams
Spectral Sequence for \boldmath$R$-modules}

\covertitle{On the Adams Spectral Sequence for $R$-modules}

\asciititle{On the Adams Spectral Sequence for R-modules}

%%%replaced
%\authors{Andrew Baker\\Andrej Lazarev}
%%%by
\authors{Andrew Baker\\Andrey Lazarev}
%%%
\address{Mathematics Department, Glasgow University,
Glasgow G12 8QW, UK. \\
Mathematics Department, Bristol University,
Bristol BS8 1TW, UK.}

\asciiemail{a.baker@maths.gla.ac.uk, a.lazarev@bris.ac.uk}
\email{a.baker@maths.gla.ac.uk\qua{\rm and}\qua  a.lazarev@bris.ac.uk}
\url{www.maths.gla.ac.uk/$\sim$ajb\qua{\rm and}\qua
www.maths.bris.ac.uk/$\sim$pure/staff/maxal/maxal}

\begin{abstract}
We discuss the Adams Spectral Sequence for $R$-modules based on commutative
localized regular quotient ring spectra over a commutative $S$-algebra $R$
in the sense of Elmendorf, Kriz, Mandell, May and Strickland. The formulation
of this spectral sequence is similar to the classical case and the calculation
of its $\E_2$-term involves the cohomology of certain `brave new Hopf algebroids'
$E^R_*E$. In working out the details we resurrect Adams' original approach to
Universal Coefficient Spectral Sequences for modules over an $R$ ring spectrum.

We show that the Adams Spectral Sequence for $S_R$ based on a commutative localized
regular quotient $R$ ring spectrum $E=R/I[X^{-1}]$ converges to the homotopy of
the $E$-nilpotent completion
\[
\pi_*\hat{\L}^R_ES_R=R_*[X^{-1}]\sphat_{I_*}.
\]
%%%replaced
%We also show that $\hat\L^R_ES_R$ is equivalent to $\L^R_ES_R$, the Bousfield
%localization of $S_R$ with respect to $E$-theory.
%%%by
We also show that when the generating regular sequence of $I_*$ is finite,
$\hat\L^R_ES_R$ is equivalent to $\L^R_ES_R$, the Bousfield localization
of $S_R$ with respect to $E$-theory.
%%%
The spectral sequence here collapses at its $\E_2$-term but it does not have
a vanishing line because of the presence of polynomial generators of positive
cohomological degree. Thus only one of Bousfield's two standard convergence
criteria applies here even though we have this equivalence. The details involve
the construction of an $I$-adic tower
\[
R/I\lla R/I^2\lla\cdots\lla R/I^s\lla R/I^{s+1}\lla\cdots
\]
whose homotopy limit is $\hat\L^R_ES_R$. We describe some examples for the
motivating case $R=\MU$.
\end{abstract}

\asciiabstract{ We discuss the Adams Spectral Sequence for R-modules
based on commutative localized regular quotient ring spectra over a
commutative S-algebra R in the sense of Elmendorf, Kriz, Mandell, May
and Strickland. The formulation of this spectral sequence is similar
to the classical case and the calculation of its E_2-term involves the
cohomology of certain `brave new Hopf algebroids' E^R_*E. In working
out the details we resurrect Adams' original approach to Universal
Coefficient Spectral Sequences for modules over an R ring spectrum.

We show that the Adams Spectral Sequence for S_R based on a
commutative localized regular quotient R ring spectrum E=R/I[X^{-1}]
converges to the homotopy of the E-nilpotent completion
pi_*hat{L}^R_ES_R=R_*[X^{-1}]^hat_{I_*}. We also show that
%%%inserted
when the generating regular sequence of I_* is finite,
%%%
hatL^R_ES_R is equivalent to L^R_ES_R, the Bousfield localization of
S_R with respect to E-theory. The spectral sequence here collapses at
its E_2-term but it does not have a vanishing line because of the presence
of polynomial generators of positive cohomological degree. Thus only one
of Bousfield's two standard convergence criteria applies here even though
we have this equivalence. The details involve the construction of an I-adic
tower
R/I <-- R/I^2 <-- ... <-- R/I^s <-- R/I^{s+1} <-- ...
whose homotopy limit is hatL^R_ES_R. We describe some examples for the
motivating case R=MU.}


\asciikeywords{S-algebra, R-module, R ring spectrum, Adams Spectral
Sequence, regular quotient}

\maketitle
\vglue12pt

\section*{Erratum}

While this paper was in e-press, the authors discovered
that the original versions of Theorems~\ref{thm:LQMU-nilcomp}
and~\ref{thm:LQ-R} were incorrect since they did not assume
that the regular sequence $u_j$ was finite. With the agreement
of the Editors, we have revised this version to include the
appropriate finiteness assumptions. We have also modified the
Abstract and Introduction to reflect this and in
Section~\ref{sec:LRQ-MU} have replaced Bousfield localizations
$\L^R_EX$ by $E$-nilpotent completions $\hat{\L}^R_EX$. As
far as we are aware, there are no further problems arising from
this mistake.
\par\vglue 5pt

\rightline{\it
Andrew Baker and Andrey Lazarev\qua 9 May 2001}

\section*{Introduction}

We consider the Adams Spectral Sequence for $R$-modules based on localized
regular quotient ring spectra over a commutative $S$-algebra $R$ in the sense
of~\cite{EKMM,Strickland:MU}, making systematic use of ideas and notation from
those two sources. This work grew out of a preprint~\cite{AB:BraveNewASS} and
the work of~\cite{AB+AJ:Brave-MU}; it is also related to ongoing collaboration
with Alain Jeanneret on Bockstein operations in cohomology theories defined on
$R$-modules~\cite{AB+AJ:Brave-BocksteinOps}.

One slightly surprising phenomenon we uncover concerns the convergence of the
Adams Spectral Sequence based on $E=R/I[X^{-1}]$, a commutative localized regular
quotient of a commutative $S$-algebra $R$. We show that the spectral sequence
for $\pi_*S_R$ collapses at $\E_2$, however for $r\geq2$, $\E_r$ has no vanishing
line because of the presence of polynomial generators of positive cohomological
degree which are infinite cycles. Thus only one of Bousfield's two convergence
criteria~\cite{Bousfield:ASS} (see Theorems~\ref{thm:BousfieldCgceASS}
and~\ref{thm:BousfieldCgceASS-Comp=Loc} below) apply here.
%%%replaced
%Despite this the
%spectral sequence converges to $\pi_*\L^R_ES_R$, where $\L^R_E$ is the
%Bousfield localization functor with respect to $E$-theory on the category
%of $R$-modules and
%%%by
Despite this, when the generating regular sequence of $I_*$ is finite,
the spectral sequence converges to $\pi_*\L^R_ES_R$, where $\L^R_E$ is
the Bousfield localization functor with respect to $E$-theory on the
category of $R$-modules and
%%%
\[
\pi_*\L^R_ES_R=R_*[X^{-1}]\sphat_{I_*},
\]
the $I_*$-adic completion of $R_*[X^{-1}]$; we also show that
%%%inserted
in this case
%%%
$\L^R_ES_R\simeq\hat\L^R_ES_R$, the $E$-nilpotent completion of $S_R$.
In the final section we describe some examples for the important case
of $R=\MU$, leaving more delicate calculations for future work.

To date there seems to have been very little attention paid to the detailed
homotopy theory associated with the category of $R$-modules, apart from general
results on Bousfield localizations and Wolbert's work on $K$-theoretic localizations
in~\cite{EKMM,Wolbert}. We hope this paper leads to further work in this area.

\subsection*{Acknowledgements}
The first author wishes to thank the University and City of Bern for providing
such a hospitable environment during many visits, also Alain Jeanneret, Urs
W\"urgler and other participants in the Topology working seminar during spring
and summer 2000. Finally, the authors wish to thank the members of Transpennine
Topology Triangle (funded by the London Mathematical Society) for providing a
timely boost to this project.

\section*{Background assumptions, terminology and technology}

We work in a setting based on a good category of spectra $\Spectra$ such
as the category of $\mathbb{L}$-spectra of~\cite{EKMM}. Associated to this
is the subcategory of $S$-modules $\mathcal M_S$ and its derived homotopy
category $\mathcal D_S$.

Throughout, $R$ will denote a commutative $S$-algebra in the sense of~\cite{EKMM}.
There is an associated subcategory $\mathcal M_R$ of $\mathcal M_S$ consisting of
the $R$-modules, and its derived homotopy category $\mathcal D_R$ and our homotopy
theoretic work is located in the latter. Because we are working in $\mathcal D_R$,
we frequently make constructions using cell $R$-modules in place of non-cell
modules (such as $R$ itself).

For $R$-modules $M$ and $N$, we set
\[
M^R_*N=\pi_*M\SSmash{R}N,\quad N_R^*M=\mathcal{D}_R(M,N)^*,
\]
where $\mathcal{D}_R(M,N)^n=\mathcal{D}_R(M,\Sigma^nN)$.

We will use the following terminology of Strickland~\cite{Strickland:MU}.
If the homotopy ring $R_*=\pi_*R$ is concentrated in even degrees, a
\emph{localized quotient} of $R$ will be an $R$ ring spectrum of the form
$R/I[X^{-1}]$. A localized quotient is \emph{commutative} if it is a commutative
$R$ ring spectrum. A localized quotient $R/I[X^{-1}]$ is \emph{regular} if the
ideal $I_*\ideal R_*$ is generated by a regular sequence $u_1,u_2,\ldots$ say.
The ideal $I_*\ideal R_*$ extends to an ideal of $R_*[X^{-1}]$ which we will
again denote by $I_*$; then as $R$-modules, $R/I[X^{-1}]\simeq R[X^{-1}]/I$.

We will make use of the language and ideas of algebraic derived categories
of modules over a commutative ring, mildly extended to deal with evenly graded
rings and their modules. In particular, this means that chain complexes are
often bigraded (or even multigraded) objects with their first grading being
homological and the second and higher ones being internal.

\section{Brave new Hopf algebroids and their cohomology}
\label{sec:Brave-HA}

If $E$ is a commutative $R$-ring spectrum, the smash product $E\SSmash{R}E$
is also a commutative $R$-ring spectrum. More precisely, it is naturally
an $E$-algebra spectrum in two ways induced from the left and right units
\[
E\xrightarrow{\;\iso\;}E\SSmash{R}R\lra E\SSmash{R}E
\lla E\SSmash{R}R\xleftarrow{\;\iso\;}E.
\]
\begin{thm}\label{thm:ERE-HA}
Let $E^R_*E$ be flat as a left or equivalently right $E_*$-module. Then
the following are true. \\
{\rm i)}
$(E_*,E^R_*E)$ is a Hopf algebroid over $R_*$. \\
{\rm ii)} for any $R$-module $M$, $E^R_*M$ is a left $E^R_*E$-comodule.
\end{thm}
\begin{proof}
This is proved using essentially the same argument as
in~\cite{Adams:Chicago,Ravenel:book}. The natural map
\[
E\SSmash{R}M\xrightarrow{\;\iso\;}E\SSmash{R}R\SSmash{R}M
\lra E\SSmash{R}E\SSmash{R}M
\]
induces the coaction
\[
\psi\:E^R_*M\lra\pi_*E\SSmash{R}E\SSmash{R}M
\xrightarrow{\;\iso\;}E^R_*E\oTimes{E_*}E^R_*M,
\]
which uses an isomorphism
\[
\pi_*E\SSmash{R}E\SSmash{R}M\iso E^R_*E\oTimes{E_*}E^R_*M.
\]
that follows from the flatness condition.
\end{proof}

For later use we record a general result on the Hopf algebroids associated
with commutative regular quotients. A number of examples for the case $R=\MU$
are discussed in Section~\ref{sec:LRQ-MU}.
\begin{prop}\label{prop:E*R/I}
Let $E=R/I$ be a commutative regular quotient where $I_*$ is generated by
the regular sequence $u_1,u_2,\ldots$. Then as an $E_*$-algebra,
\[
E^R_*E=\Lambda_{E_*}(\tau_i:i\geq1),
\]
where $\deg\tau_i=\deg u_i+1$. Moreover, the generators $\tau_i$ are
primitive with respect to the coaction, and $E^R_*E$ is a primitively
generated Hopf algebra over $E_*$.

Dually, as an $E_*$-algebra,
\[
E_R^*E=\hat\Lambda_{E_*}(Q^i:i\geq1),
\]
where $Q^i$ is the Bockstein operation dual to $\tau_i$ with
$\deg Q^i=\deg u_i+1$ and $\hat\Lambda_{E_*}(\ )$ indicates the
completed exterior algebra generated by the anti-commuting $Q^i$
elements.
\end{prop}

The proof requires the K\"unneth Spectral Sequence for $R$-modules
of~\cite{EKMM},
\[
\E^2_{p,q}=\Tor^{R_*}_{p,q}(E_*,E_*)\Lra E^R_{p+q}E.
\]
This spectral sequence is multiplicative, however there seems to be no
published proof in the literature. At the suggestion of the referee, we
indicate a proof of this due to M.~Mandell and which originally appeared
in a preprint version of~\cite{Lazarev}.
\begin{lem}\label{lem:KSS-Mult}
If $A$ and $B$ are $R$ ring spectra then the K\"unneth Spectral Sequence
\[
\Tor^{R_*}(A_*,B_*)\Lra A^R_*B=\pi_*A\SSmash{R}B
\]
is a spectral sequence of differential graded $R_*$-algebras.
\end{lem}
\begin{proof}[Sketch proof]
To deal with the multiplicative structure we need to modify the
original construction given in Part~IV section~5 of~\cite{EKMM}.
We remind the reader that we are working in the derived homotopy
category $\mathcal D_R$.

Let
\[
\cdots\lra F_{p,*}\xrightarrow{f_p}F_{p-1,*}\lra\cdots
\xrightarrow{f_1}F_{0,*}\xrightarrow{f_0}A_*\ra0
\]
be an free $R_*$-resolution of $A_*$. Using freeness, we can choose
a map of complexes
\[
\mu\:F_{*,*}\oTimes{R_*}F_{*,*}\lra F_{*,*}
\]
which lifts the multiplication on $A_*$.

For each $p\geq 0$ let $\mathbf F_p$ be a wedge of sphere $R$-modules
satisfying $\pi_*\mathbf F_p=F_{p,*}$. Set $A'_0=\mathbf F_0$ and choose
a map $\phi_0\:A'_0\lra A$ inducing $f_0$ in homotopy. If $\mathbf Q_0$
is the homotopy fibre of $\phi_0$ then
\[
\pi_*\mathbf Q_0=\ker f_0
\]
and we can choose a map $\mathbf F_1\lra\mathbf Q_0$ for which the composition
$\phi'_1\:\mathbf F_1\lra\mathbf Q_0\lra\mathbf F_0$ induces $f_1$ in homotopy.
Next take $A'_1$ to be the cofibre of $\phi'_1$. The map $\phi_0$ has a canonical
extension to a map $\phi_1\:A'_1\lra A$. If $\mathbf Q_1$ is the homotopy fibre
of $\phi_1$ then
\[
\pi_*\Sigma^{-1}\mathbf Q_1=\ker f_1,
\]
and we can find a map $\mathbf F_2\lra\Sigma\mathbf Q_1$ for which the composite
map $\phi'_2\:\mathbf F_2\lra\mathbf Q_1\lra\mathbf F_1$ induces $f_2$ in homotopy.
We take $A'_2$ to be the cofibre of $\phi'_2$ and find that there is a canonical
extension of $\phi_1$ to a map $\phi_2\:A'_2\lra A$.

Continuing in this way we construct a directed system
\begin{equation}\label{eqn:fil}
A'_0\lra A'_1\lra\cdots\lra A'_p\lra\cdots
\end{equation}
whose telescope $A'$ is equivalent to $A$. Since we can assume that all consecutive
maps are inclusions of cell subcomplexes, there is an associated filtration on $A'$.
Smashing this with $B$ we get a filtration on $A'\SSmash{R}B$ and an associated
spectral sequence converging to $A_*^RB$. The identification of the $\E_2$-term is
routine.

Recall that $A$ and therefore $A'$ are $R$ ring spectra. Smashing the directed
system of~\eqref{eqn:fil} with itself we obtain a filtration on $A'\SSmash{R}A'$,
\begin{equation}\label{eqn:mul}
A'_0\SSmash{R}A'_0\lra\cdots\lra\bigcup_{i+j=k}A'_i\SSmash{R}A'_j
\lra\bigcup_{i+j=k+1}A'_i\SSmash{R}A'_j\lra\cdots,
\end{equation}
where the filtrations terms are unions of the subspectra $A_i'\SSmash{R}A_j'$.
Proceeding by induction, we can realize the multiplication map $A'\SSmash{R}A'\lra A'$
as a map of filtered $R$-modules so that on the cofibres of the filtration terms
of~\eqref{eqn:mul} it agrees with the pairing $\mu$.

We have constructed a collection of maps $A'_i\SSmash{R}A'_j\lra A'_{i+j}$.
Using these maps and the multiplication on $B$ we can now construct maps
\[
A'_i\SSmash{R}B\SSmash{R}A'_j\SSmash{R}B\lra A_{i+j}\SSmash{R}B
\]
which induce the required pairing of spectral sequences.
\end{proof}
\begin{proof}[Proof of Proposition~\ref{prop:E*R/I}]
As in the discussion preceding Proposition~\ref{prop:Tor(R/I,R/I)},
making use of a Koszul resolution we obtain
\[
\E^2_{*,*}=\Lambda_{E_*}(e_i:i\geq1).
\]
The generators have bidegree $\bideg e_i=(1,|u_i|)$, so the differentials
\[
\d^r\:\E^r_{p,q}\lra\E^r_{p-r,q+r-1}
\]
are trivial on the generators $e_i$ for dimensional reasons. Together with
multiplicativity, this shows that spectral sequence collapses, giving
\[
E^R_*E=\Lambda_{E_*}(\tau_i:i\geq1),
\]
where the generator $\tau_i$ has degree $\deg\tau_i=\deg u_i+1$ and is
represented by $e_i$.

For each $i$,
\[
(R/u_i)^R_*(R/u_i)=\Lambda_{R_*/(u_i)}(\tau'_i)
\]
with $\deg\tau'_i=|u_i|+1$. Under the coproduct, $\tau'_i$ is primitive
for degree reasons. By comparing the two K\"unneth Spectral Sequences we
find that $\tau_i\in E^R_*E$ can be chosen to be the image of $\tau'_i$
under the evident ring homomorphism $(R/u_i)^R_*(R/u_i)\lra E^R_*E$, which
is actually a morphism of Hopf algebroids over $R_*$. Hence $\tau_i$ is
coaction primitive in $E^R_*E$.

For $E_R^*E$, we construct the Bockstein operation $Q^i$ using the composition
\[
R/u_i\lra\Sigma^{|u_i|+1}R\lra\Sigma^{|u_i|+1}R/u_i
\]
to induce a map $E\lra\Sigma^{|u_i|+1}E$, then use the Koszul resolution
to determine the Universal Coefficient Spectral sequence
\[
\E_2^{p,q}=\Ext_{R_*}^{p,q}(E_*,E_*)\Lra E_R^{p+q}E
\]
which collapses at its $\E_2$-term. Further details on the construction
of these operations appear in~\cite{Strickland:MU,AB+AJ:Brave-BocksteinOps}.
\end{proof}
\begin{cor}\label{cor:E*R/I}
{\rm i)} The natural map $E_*=E^R_*R\lra E^R_*E$ induced by the unit
$R\lra R/I$ is a split monomorphism of $E_*$-modules. \\
{\rm ii)} $E^R_*E$ is a free $E_*$-module.
\end{cor}
\begin{proof}
An explicit splitting as in (i) is obtained using the multiplication
map $E\SSmash{R}E\lra E$ which induces a homomorphism of $E_*$-modules
$E^R_*E\lra E_*$.
\end{proof}

We will use $\Coext$ to denote the cohomology of such Hopf algebroids
rather than $\Ext$ since we will also make heavy use of $\Ext$ groups
for modules over rings; more details of the definition and calculations
can be found in~\cite{Adams:Chicago,Ravenel:book}. Recall that for
$E^R_*E$-comodules $L_*$ and $M_*$ where $L_*$ is $E_*$-projective,
$\Coext_{E^R_*E}^{s,t}(L_*,M_*)$ can be calculated as follows. Consider
a resolution
\[
0\ra M_*\lra J_{0,*}\lra J_{1,*}\lra\cdots\lra J_{s,*}\lra\cdots
\]
in which each $J_{s,*}$ is a summand of an extended comodule
\[
E^R_*E\Square{E_*}N_{s,*},
\]
for some $E_*$-module $N_{s,*}$. Then the complex
\begin{multline*}
0\ra\Hom_{E^R_*E}^*(L_*,J_{0,*})\lra\Hom_{E^R_*E}^*(L_*,J_{1,*})
\\
\lra\cdots\lra\Hom_{E^R_*E}^*(L_*,J_{s,*})
\lra\cdots
\end{multline*}
has cohomology
\[
\mathrm{H}^s(\Hom_{E^R_*E}^*(L_*,J_{*,*}))=\Coext_{E^R_*E}^{s,*}(L_*,M_*).
\]
The functors $\Coext_{E^R_*E}^{s,*}(L_*,\ )$ are the right derived functors
of the left exact functor
\[
M_*\rightsquigarrow\Hom_{E^R_*E}^*(L_*,M_*)
\]
on the category of left $E^R_*E$-comodules. By analogy with~\cite{Ravenel:book},
when $L_*=E_*$ we have
\[
\Coext_{E^R_*E}^{s,*}(E_*,M_*)=\Cotor_{E^R_*E}^{s,*}(E_*,M_*).
\]

\section{The Adams Spectral Sequence for $R$-modules}\label{sec:ASS}

We will describe the $E$-theory Adams Spectral Sequence in the homotopy
category of $R$-module spectra. As in the classical case of sphere spectrum
$R=S$, it turns out that the $\E_2$-term is can be described in terms of the
functor $\Coext_{E^R_*E}$.

Let $L,M$ be $R$-modules and $E$ a commutative $R$-ring spectrum with
$E^R_*E$ flat as a left (or right) $E_*$-module.
\begin{thm}\label{thm:ASS}
If $E^R_*L$ is projective as an $E_*$-module, there is an
Adams Spectral Sequence with
\[
\E_2^{s,t}(L,M)=\Coext_{E^R_*E}^{s,t}(E^R_*L,E^R_*M).
\]
\end{thm}
\begin{proof}
Working throughout in the derived category $\mathcal{D}_R$, the proof
follows that of Adams~\cite{Adams:Chicago}, with $S_R\simeq R$ replacing
the sphere spectrum $S$. The canonical Adams resolution of $M$ is built
up in the usual way by splicing together the cofibre triangles in the
following diagram.
\[
\begin{xy}
\xymatrix@-0.7pc{
\ar[dr]M&&\bar{E}\SSmash{R}\ar[ll]M
\ar[dr]&&\bar{E}\SSmash{R}\bar{E}\SSmash{R}\ar[ll]M
\ar[dr]&\ar[l]&\ar@{.}[l]
&&
\\
&E\SSmash{R}M\ar[ur]&
&E\SSmash{R}\bar{E}\SSmash{R}M\ar[ur]&
&&
&&
}
\end{xy}
\]
The algebraic identification of the $\E_2$-term proceeds as
in~\cite{Adams:Chicago}.
\end{proof}
In the rest of this paper we will have $L=S_R\simeq R$, and set
\[
\E_2^{s,t}(M)=\Coext_{E^R_*E}^{s,t}(E_*,E^R_*M).
\]
We will refer to this spectral sequence as the Adams Spectral
Sequence based on $E$ for the $R$-module $M$.

To understand convergence of such a spectral sequence we use
a criterion of Bousfield~\cite{Bousfield:ASS,Ravenel:LocnPaper}.
For an $R$-module $M$, let $D_sM$ ($s\geq0$) be the $R$-modules
defined by $D_0M=M$ and taking $D_sM$ to be the fibre of the
natural map
\[
D_{s-1}M\iso R\SSmash{R}D_{s-1}M\lra E\SSmash{R}D_{s-1}M.
\]
Also for each $s\geq0$ let $K_sM$ be the cofibre of the natural map
$D_sM\lra M$. Then the \emph{$E$-nilpotent completion} of $M$ is the
homotopy limit
\[
\hat{\L}^R_EM=\holim_sK_sM.
\]
\begin{rem}\label{rem:Ehat-Invce}
It is easy to see that if $M\lra N$ is a map of $R$-modules which
is an $E$-equivalence, then for each $s$, there is an equivalence
$K_sM\lra K_sN$, hence
\[
\hat{\mathrm{L}}^R_EM\simeq\hat{\mathrm{L}}^R_EN.
\]
\end{rem}
\begin{thm}\label{thm:BousfieldCgceASS}
If for each pair $(s,t)$ there is an $r_0$ for which
$\E_r^{s,t}(M)=\E_\infty^{s,t}(M)$ whenever $r\geq r_0$,
then the Adams Spectral Sequence for $M$ based on $E$
converges to $\pi_*\hat{\mathrm{L}}^R_EM$.
\end{thm}
Although there is a natural map $\mathrm{L}^R_EM\lra\hat{\mathrm{L}}^R_EM$,
it is not in general a weak equivalence; this equivalence is guaranteed
by another result of Bousfield~\cite{Bousfield:ASS}.
\begin{thm}\label{thm:BousfieldCgceASS-Comp=Loc}
Suppose that there is an $r_1$ such that for every $R$-module $N$ there
is an $s_1$ for which $\E_{r}^{s,t}(N)=0$ whenever $r\geq r_1$ and
$s\geq s_1$. Then for every $R$-module $M$ the Adams Spectral Sequence
for $M$ based on $E$ converges to $\pi_*\mathrm{L}^R_EM$ and
\[
\mathrm{L}^R_EM\simeq\hat{\mathrm{L}}^R_EM.
\]
\end{thm}

\section{The Universal Coefficient Spectral Sequence for regular quotients}
\label{sec:UCKSS-LQR}

Let $R$ be a commutative $S$-algebra and $E=R/I$ a commutative regular
quotient of $R$, where $u_1,u_2,\ldots$ is a regular sequence generating
$I_*\ideal R_*$.

We will discuss the existence of the Universal Coefficient Spectral Sequence
\begin{equation}\label{eqn:E*-UCSS}
\E^2_{r,s}=\Ext^{r,s}_{E_*}(E^R_*M,N_*)\Lra N_R^*M,
\end{equation}
where $M$ and $N$ are $R$-modules and $N$ is also an $E$-module spectrum in
$\mathcal M_R$. The classical prototype of this was described by
Adams~\cite{Adams:Chicago} (who generalized a construction of Atiyah~\cite{Atiyah}
for the K\"unneth Theorem in $K$-theory) and used in setting up the $E$-theory
Adams Spectral Sequence. It is routine to verify that Adams' approach can be
followed in $\mathcal D_R$. We remark that if $E$ were a commutative $R$-algebra
then the Universal Coefficient Spectral Sequence of~\cite{EKMM} would be applicable
but that condition does not hold in the generality we require.

The existence of such a spectral sequence depends on the following
conditions
being satisfied.
\begin{cond}\label{cond:UCSS}
$E$ is a homotopy colimit of finite cell $R$-modules $E_\alpha$ whose
$R$-Spanier Whitehead duals $\mathrm{D}_RE_\alpha=\Func_R(E_\alpha,R)$
satisfy the two conditions  \\
(A) $E^R_*\mathrm{D}_RE_\alpha$ is $E_*$-projective; \\
(B) the natural map
\[
N_R^*M\lra\Hom_{E_*}(E^R_*M,N_*)
\]
is an isomorphism.
\end{cond}
\begin{thm}\label{thm:UCSS-RQ-Conditions}
For a commutative regular quotient $E=R/I$ of $R$, $E$ can be expressed as
a homotopy colimit of finite cell $R$-modules satisfying the conditions of
\emph{Condition~\ref{cond:UCSS}}. In fact we can take
$E^R_*\mathrm{D}_RE_\alpha$ to be $E_*$-free.
\end{thm}
The proof will use the following Lemma.
\begin{lem}\label{lem:UCSS-lemma}
Let $u\in R_{2d}$ be non-zero divisor in $R_*$. Suppose that $P$ is an
$R$-module for which $E^R_*P$ is $E_*$-projective and for an $E$-module
$R$-spectrum $N$,
\[
N_R^*P\iso\Hom_{E_*}(E^R_*P,N_*).
\]
Then $E^R_*P\SSmash{R}R/u$ is $E_*$-projective and
\[
N_R^*P\SSmash{R}R/u\iso\Hom_{E_*}(E^R_*P\SSmash{R}R/u,N_*).
\]
\end{lem}
\begin{proof}
Smashing $E\SSmash{R}P$ with the cofibre sequence~\eqref{eqn:UCSS-CofibSeq}
and taking homotopy, we obtain an exact triangle
\[
\begin{xy}
\xymatrix{
E^R_*P\ar[rr]^{\ds u}&&E^R_*P\ar[dl]\\
&E^R_*P\SSmash{R}R/u\ar[ul]&
}
\end{xy}
\]
As multiplication by $u$ induces the trivial map in $E^R$-homology,
this is actually a short exact sequence of $E_*$-modules,
\[
0\ra E^R_*P\lra E^R_*P\SSmash{R}R/u\lra E^R_*P\ra0
\]
which clearly splits, so $E^R_*P\SSmash{R}R/u$ is $E_*$-projective.

In the evident diagram of exact triangles
\[
\begin{xy}
\xymatrix@-1pc{
N_R^*P\ar[rr]\ar[dd]&&N_R^*P\ar[dl]\ar[dd]   \\
&N_R^*P\SSmash{R}R/u\ar[ul]\ar[dd]&   \\
\Hom_{E_*}(E^R_*P,N_*)\ar'[r][rr]&&\Hom_{E_*}(E^R_*P,N_*)\ar[dl] \\
&\Hom_{E_*}(E^R_*P\SSmash{R}R/u,N_*)\ar[ul]&
}
\end{xy}
\]
the map $N_R^*P\lra\Hom_{E_*}(E^R_*P,N_*)$ is an isomorphism, so
\[
N_R^*P\SSmash{R}R/u\lra\Hom_{E_*}(E^R_*P\SSmash{R}R/u,N_*)
\]
is also an isomorphism by the Five Lemma.
\end{proof}
\begin{proof}[Proof of Theorem~\ref{thm:UCSS-RQ-Conditions}]
Let $u_1,u_2,\ldots$ be a regular sequence generating $I_*\ideal R_*$.
Using the notation $R/u=R/(u)$, we recall from~\cite{Strickland:MU}
that
\[
E=\hocolim_k R/u_1\SSmash{R}R/u_2\SSmash{R}\cdots\SSmash{R}R/u_k.
\]
For $u\in R_{2d}$ a non-zero divisor, the $R_*$-free resolution
\[
0\ra R_*\xrightarrow{\ph{\;u\;}}R_*\xrightarrow{\;u\;}R_*/(u)\ra0
\]
corresponds to an $R$-cell structure on $R/u$ with one cell in each of
the dimensions $0$ and $2d+1$. There is an associated cofibre sequence
\begin{equation}\label{eqn:UCSS-CofibSeq}
\cdots\lra
\Sigma^{2d}R\xrightarrow{u}R\lra R/u\lra\Sigma^{2d+1}R
\lra\cdots,
\end{equation}
for which the induced long exact sequence in $E^R$-homology shows that
$E^R_*R/u$ is $E_*$-free. The dual $\mathrm{D}_RR/u$ is equivalent to
$\Sigma^{-(2d+1)}R/u$, hence $R/u$ is essentially self dual.

For an $E$-module spectrum $N$ in $\mathcal D_R$, there are two exact
triangles and morphisms between them,
\[
\begin{xy}
\xymatrix@-1pc{
N_R^*R\ar[rr]\ar[dd]&&N_R^*R\ar[dl]\ar[dd] \\
&N_R^*R/u\ar[ul]\ar[dd]&  \\
\Hom_{E_*}(E_*,N_*)\ar'[r][rr]&&\Hom_{E_*}(E_*,N_*)\ar[dl] \\
& \Hom_{E_*}(E^R_*R/u,N_*)\ar[ul]&
}
\end{xy}
\]
The identifications
\[
N_*\iso N_R^*R\iso\Hom_{E_*}(E_*,N_*),
\]
and the Five Lemma imply that
\[
N_R^*R/u\iso\Hom_{E_*}(E^R_*R/u,N_*).
\]
Lemma~\ref{lem:UCSS-lemma} now implies that each of the spectra
$R/u_1\SSmash{R}R/u_2\SSmash{R}\cdots\SSmash{R}R/u_k$ satisfies
conditions (A) and (B).
\end{proof}

\section{The Adams Spectral Sequence based on a regular quotient}
\label{sec:ASS-RQ}

For an $R$-module $M$, let $M^{(s)}$ denote the $s$-fold $R$-smash
power of $M$,
\[
M^{(s)}=M\SSmash{R}M\SSmash{R}\cdots\SSmash{R}M.
\]
If $M$ is an $R[X^{-1}]$-module, then
\[
M^{(s)}=
M\SSmash{R[X^{-1}]}M\SSmash{R[X^{-1}]}\cdots\SSmash{R[X^{-1}]}M.
\]

Let $E=R/I[X^{-1}]$ be a localized regular quotient and $u_1,u_2,\ldots$
a regular sequence generating $I_*$. We will discuss the Adams Spectral
Sequence based on $E$. By Remark~\ref{rem:Ehat-Invce}, we can work in the
category of $R[X^{-1}]$-modules and replace the Adams Spectral Sequence of
$S_R$ by that of $S_{R[X^{-1}]}$. To simplify notation, from now on we will
replace $R$ by $R[X^{-1}]$ and therefore assume that $E=R/I$ is a regular
quotient of $R$.

First we identify the canonical Adams resolution giving rise to the Adams
Spectral Sequence based on the regular quotient $E=R/I$. We will relate this
to a tower described by the second author~\cite{Lazarev}, but the reader
should beware that his notation for $I^{(s)}$ is $I^s$ which we will use for
a different spectrum.

There is a fibre sequence $I\lra R\lra R/I$ and a tower of maps of $R$-modules
\[
R\lla I\lla I^{(2)}\lla\cdots\lla I^{(s)}\lla I^{(s+1)}\lla\cdots
\]
in which $I^{(s+1)}\lra I^{(s)}$ is the evident composite
\[
I^{(s+1)}\lra R\SSmash{R}I^{(s)}=I^{(s)}.
\]
Setting $R/I^{(s)}=\cofibre(I^{(s)}\lra R)$, we obtain a tower
\[
R/I\lla R/I^{(2)}\lla\cdots\lla R/I^{(s)}\lla R/I^{(s+1)}\lla\cdots
\]
which we will refer to as the \emph{external $I$-adic tower}. The
next result is immediate from the definitions.
\begin{prop}\label{prop:ASS-LRQ-Tower}
We have
\[
D_0S_{R}=R,
\qquad
D_sS_{R}=I^{(s)},
\quad(s\geq1),
\]
and
\[
K_sS_R=R/I^{(s+1)}
\quad(s\geq0).
\]
\end{prop}

It is not immediately clear how to determine the limit
\[
\hat{\mathrm{L}}^R_ES_R=\holim_{s}R/I^{(s)}.
\]
Instead of doing this directly, we will adopt an approach suggested by
Bousfield~\cite{Bousfield:ASS}, making use of another $E$-nilpotent
resolution, associated with the \emph{internal $I$-adic tower} to be
described below.

In order to carry this out, we first need to understand convergence.
We will see that the condition of Theorem~\ref{thm:BousfieldCgceASS}
is satisfied for a commutative regular quotient $E=R/I$.
\begin{prop}\label{prop:LRQ-ASS-E2}
The $\E_2$-term of the $E$-theory Adams Spectral Sequence for $\pi_*S_R$
is
\[
\E_2^{s,t}(S_R)=\Coext_{E^R_*E}^{s,t}(E_*,E_*)
=E_*[U_i:i\geq1],
\]
where $\bideg U_i=(1,|u_i|+1)$. Hence this spectral sequence
collapses at its $\E_2$-term
\[
\E_2^{*,*}(S_R)=\E_\infty^{*,*}(S_R)
\]
and converges to $\pi_*\hat{\mathrm{L}}^R_ES_R$.
\end{prop}
\begin{proof}
By Proposition~\ref{prop:E*R/I},
\[
E^R_*E=\Lambda_{R_*}(\tau_i:i\geq1),
\]
with generators $\tau_i$ which are primitive with respect to the
coproduct of this Hopf algebroid. The determination of
\[
\Coext_{E^R_*E}^{*,*}(E_*,E_*)
\]
is now standard and the differentials are trivial for degree reasons.
\end{proof}
Induction on the number of cells now gives
\begin{cor}\label{cor:LRQ-ASS-E2}
For a finite cell $R$-module $M$, the $E$-theory Adams Spectral Sequence
for $\pi_*M$ converges to $\pi_*\hat{\L}^R_EM$.
\end{cor}

\section{The internal $I$-adic tower}\label{sec:I-adictower}

Suppose that $I_*\ideal R_*$ is generated by a regular sequence
$u_1,u_2,\ldots$. We will often indicate a monomial in the $u_i$ by writing
$u_{(i_1,\ldots,i_k)}=u_{i_1}\cdots u_{i_k}$. We will write $E=R/I$ and make
use of algebraic results from~\cite{AB:HomRegQuot} which we now recall in detail.

For $s\geq0$, we define the $R$-module $I^s/I^{s+1}$ to be the wedge of copies
of $E$ indexed on the distinct monomials of degree $s$ in the generators $u_i$.
For an explanation of this, see Corollary~\ref{cor:I*s/I*(s+1)-R/Ifree}.

We will show that there is an (\emph{internal}) \emph{$I$-adic tower}
of $R$-modules
\[
R/I\lla R/I^2\lla\cdots\lla R/I^s\lla R/I^{s+1}\lla\cdots
\]
so that for each $s\geq0$ the fibre sequence
\[
R/I^s\lla R/I^{s+1}\lla I^s/I^{s+1}
\]
corresponds to a certain element of
\[
\Ext_{R_*}^{1}(R_*/I_*^s,I_*^s/I_*^{s+1})
\]
in $\E_2$-term of the Universal Coefficient Spectral Sequence
of~\cite{EKMM} converging to $\mathcal{D}_R(R/I^s,I^s/I^{s+1})^*$.
On setting $I^s=\fibre(R\lra R/I^s)$ we obtain another tower
\[
R\lla I\lla I^2\lla\cdots\lla I^s\lla I^{s+1}\lla\cdots
\]
which is analogous to the external version of~\cite{Lazarev}.
A related construction appeared in~\cite{AB:Ainfty,AB-UW:Bockop}
for the case of $R=\hat{E(n)}$ (which was shown to admit a not
necessarily commutative $S$-algebra structure) and $I=I_n$.

Underlying our work is the classical \emph{Koszul resolution}
\[
\mathbf{K}_{*,*}\lra R_*/I_*\ra0,
\]
where
\[
\mathbf{K}_{*,*}=\Lambda_{R_*}(e_i:i\geq1),
\]
which has grading given by $\deg e_i=|u_i|+1$ and differential
\begin{align*}
\d e_i&=u_i, \\
\d(xy)&=(\d x)y+(-1)^{r}x\d y
\quad(x\in\mathbf{K}_{r,*},\;y\in\mathbf{K}_{s,*}).
\end{align*}
Hence $(\mathbf{K}_{*,*},\d)$ is an $R_*$-free resolution of $R_*/I_*$
which is a differential graded $R_*$-algebra. Tensoring with $R_*/I_*$
and taking homology leads to a well known result.
\begin{prop}\label{prop:Tor(R/I,R/I)}
As an $R_*/I_*$-algebra,
\[
\Tor^{R_*}_{*,*}(R_*/I_*,R_*/I_*)=\Lambda_{R_*/I_*}(e_i:i\geq1).
\]
\end{prop}
\begin{cor}\label{cor:Tor(R/I,R/I)}
$\Tor^{R_*}_{*,*}(R_*/I_*,R_*/I_*)$ is a free $R_*/I_*$-module.
\end{cor}
This is of course closely related to the topological result
Proposition~\ref{prop:E*R/I}.

Now returning to our algebraic discussion, we recall the following
standard result.
\begin{lem}[\cite{Matsumura}, Theorem 16.2]\label{lem:I*s/I*(s+1)-R/Ifree}
For $s\geq0$, $I_*^s/I_*^{s+1}$ is a free $R_*/I_*$-module with a
basis consisting of residue classes of the distinct monomials
$u_{(i_1,\ldots,i_s)}$ of degree $s$.
\end{lem}
\begin{cor}\label{cor:I*s/I*(s+1)-R/Ifree}
For $s\geq0$, there is an isomorphism of $R_*$-modules
\[
\pi_*I^s/I^{s+1}=I_*^s/I_*^{s+1}.
\]
Hence $\pi_*I^s/I^{s+1}$ is a free $R_*/I_*$-module with a basis indexed
on the distinct monomials $u_{(i_1,\ldots,i_s)}$ of degree $s$.
\end{cor}

Let $\mathrm{U}^{(s)}_*$ be the free $R_*$-module on a basis indexed on
the distinct monomials of degree $s$ in the $u_i$. For $s\geq0$, set
\[
\mathbf{Q}^{(s)}_{*,*}=\mathbf{K}_{*,*}\oTimes{R_*}\mathrm{U}^{(s)}_*,
\quad
\d_{\mathbf{Q}}^{(s)}=\d\otimes1,
\]
and also for $x\in\mathbf{K}_{*,*}$ write
\[
x\tilde u_{(i_1,\ldots,i_s)}=x\otimes u_{(i_1,\ldots,i_s)}.
\]
There is an obvious augmentation
\[
\mathbf{Q}^{(s)}_{0,*}\lra I_*^s/I_*^{s+1}.
\]
\begin{lem}\label{lem:I*s/I*(s+1)}
For $s\geq1$,
\[
\mathbf{Q}^{(s)}_{*,*}\xrightarrow{\epsilon^{(s)}}I_*^s/I_*^{s+1}\ra0
\]
is a resolution by free $R_*$-modules.
\end{lem}

Given a complex $(\mathbf{C}_{*,*},\d_{\mathbf{C}})$, the $k$-shifted
complex $(\mathbf{C}[k]_{*,*},\d_{\mathbf{C}[k]})$ is defined by
\[
\mathbf{C}[k]_{n,*}=\mathbf{C}_{n+k,*},
\quad
\d_{\mathbf{C}[k]}=(-1)^k\d_{\mathbf{C}}.
\]

There is a morphism of chain complexes
\begin{align*}
\del^{(s+1)}\:&\mathbf{Q}^{(s)}_{*,*}\lra\mathbf{Q}^{(s+1)}[-1]_{*,*}; \\
\del^{(s+1)}e_{i_1}\cdots e_{i_r}\tilde u_{(j_1,\ldots,j_s)}&=
\sum_{k=1}^r(-1)^k
e_{i_1}\cdots\hat e_{i_k}\cdots e_{i_r}\tilde u_{(j_1,\ldots,j_s)}.
\end{align*}
Using the identification
$\mathbf{Q}^{(s+1)}[-1]_{n,*}=\mathbf{Q}^{(s+1)}_{n-1,*}$,
we will often view $\del^{(s+1)}$ as a homomorphism
\[
\del^{(s+1)}\:\mathbf{Q}^{(s)}_{*,*}\lra\mathbf{Q}^{(s+1)}_{*,*}
\]
of bigraded $R_*$-modules of degree $-1$.

There are also external pairings
\begin{align*}
\mathbf{Q}^{(r)}_{*,*}\oTimes{R_*}\mathbf{Q}^{(s)}_{*,*}
& \lra\mathbf{Q}^{(r+s)}_{*,*}; \\
x\tilde u_{(i_1,\ldots,i_s)}\otimes y\tilde u_{(j_1,\ldots,j_s)}
&\longmapsto
xy\tilde u_{(i_1,\ldots,i_s,j_1,\ldots,j_s)}
\quad(x,y\in\mathbf{K}_{*,*}).
\end{align*}
In particular, each $\mathbf{Q}^{(r)}_{*,*}$ is a differential module
over the differential graded $R_*$-algebra $\mathbf{K}^{(0)}_{*,*}$
and $\del^{(s+1)}$ is a $\mathbf{K}^{(0)}_{*,*}$-derivation.
\begin{thm}\label{thm:R/Is-Resn}
For $s\geq1$, there is a resolution
\[
\mathbf{K}^{(s-1)}_{*,*}\xrightarrow{\epsilon^{(s-1)}}R_*/I_*^s\ra0,
\]
by free $R_*$-modules, where
\[
\mathbf{K}^{(s-1)}_{*,*}=
\mathbf{Q}^{(0)}_{*,*}\oplus\mathbf{Q}^{(1)}_{*,*}
\oplus\cdots\oplus\mathbf{Q}^{(s-1)}_{*,*},
\]
and the differential is
\[
\d^{(s-1)}=
(\d_{\mathbf{Q}}^{(0)},\del^{(1)}+\d_{\mathbf{Q}}^{(1)},
\del^{(2)}+\d_{\mathbf{Q}}^{(2)},\ldots,
\del^{(s-1)}+\d_{\mathbf{Q}}^{(s-1)}).
\]
In fact $(\mathbf{K}^{(s-1)}_{*,*},\d^{(s-1)})$ is a differential
graded $R_*$-algebra which provides a multiplicative resolution
of $R_*/I^s$, with the augmentation given by
\[
\epsilon^{(s-1)}
(x_0,x_1\tilde u_{\mathbf{i}_1},\ldots,x_{s-1}\tilde u_{\mathbf{i}_{s-1}})
=
x_0+x_1u_{\mathbf{i}_1}+\cdots+x_{s-1}u_{\mathbf{i}_{s-1}}.
\]
\end{thm}

The algebraic extension of $R_*$-modules
\[
0\la R_*/I_*^s\lla R_*/I_*^{s+1}\lla I_*^s/I_*^{s+1}\la0
\]
is classified by an element of
\[
\Ext_{R_*}^{1}(R_*/I_*^s,I_*^s/I_*^{s+1})=
\Hom_{\mathcal{D}_{R_*}}(R_*/I_*^s,I_*^s/I_*^{s+1}[-1]),
\]
where $\Hom_{\mathcal{D}_{R_*}}$ denotes morphisms in the derived category
$\mathcal{D}_{R_*}$ of the ring $R_*$~\cite{Weibel}. This element is
represented by the composite
\begin{equation}\label{eqn:R/Is-R/I(s+1)-Ext1}
\tilde\del^{(s)}_*\:
\mathbf{K}^{(s-1)}_{*,*}\xrightarrow{\text{proj}}
\mathbf{Q}^{(s-1)}_{*,*}\xrightarrow{\del^{(s)}}\mathbf{Q}^{(s)}[-1]_{*,*}.
\end{equation}

The analogue of the next result for ungraded rings was proved
in~\cite{AB:HomRegQuot}; the proof is easily adapted to the
graded case.
\begin{prop}\label{prop:QComp-Exactness}
For each $s\geq2$, the following complex is exact:
\begin{multline*}
\Tor^{R_*}_{*,*}(R_*/I_*,R_*/I_*)\xrightarrow{\del^{(1)}_*}
\Tor^{R_*}_{*,*}(R_*/I_*,I_*/I_*^2) \\
\xrightarrow{\del^{(2)}_*}
\cdots\xrightarrow{\del^{(s-1)}_*}
\Tor^{R_*}_{*,*}(R_*/I_*,I_*^{s-1}/I_*^{s}).
\end{multline*}
\end{prop}
\begin{thm}\label{thm:Tor(R/I,R/Is)}
For $s\geq2$,
\[
\Tor^{R_*}_{*,*}(R_*/I_*,R_*/I_*^s)=R_*/I_*\oplus\coker\del^{(s-1)}_*.
\]
This is a free $R_*/I_*$-module and with its natural $R_*/I_*$-algebra
structure, $\Tor^{R_*}_{*,*}(R_*/I_*,R_*/I_*^s)$ has trivial products.
\end{thm}
Given this algebraic background, we can now construct the $I$-adic
tower.
\begin{thm}\label{thm:R/I^s}
There is a tower of $R$-modules
\[
R/I\lla R/I^2\lla\cdots\lla R/I^s\lla R/I^{s+1}\lla\cdots
\]
whose maps define fibre sequences
\[
R/I^s\lla R/I^{s+1}\lla I^s/I^{s+1}
\]
which in homotopy realise the exact sequences of $R_*$-modules
\[
0\la R_*/I_*^s\lla R_*/I_*^{s+1}\lla I_*^s/I_*^{s+1}\la0.
\]
Furthermore, the following conditions are satisfied for each $s\geq1$.\\
{\rm(i)} $E^R_*R/I^s$ is a free $E_*$-module and the unit induces
a splitting
\[
E^R_*R/I^s=E_*\oplus(\ker\:E^R_*R/I^s\lra E_*);
\]
{\rm(ii)} the projection map $R/I^{s+1}\lra R/I^s$ induces
the zero map
\[
(\ker\:E^R_*R/I^{s+1}\lra E_*)\lra(\ker\:E^R_*R/I^s\lra E_*);
\]
{\rm(iii)} the inclusion map $j_s\:I^s/I^{s+1}\lra R/I^{s+1}$
induces an exact sequence
\[
E^R_*I^{s-1}/I^s\xrightarrow{\del^{(s)}_*}E^R_*I^s/I^{s+1}
\xrightarrow{{j_s}_*}(\ker\:E^R_*R/I^{s+1}\lra E_*)\ra0.
\]
\end{thm}
\begin{proof}
The proof is by induction on $s$. Assuming that $R/I^s$ exists
with the asserted properties, we will define a suitable map
$\delta_s\:R/I^s\lra\Sigma I^s/I^{s+1}$ which induces a fibre
sequence of the form
\begin{equation}\label{eqn:FibSeq(s+1)}
R/I^s\lla X^{(s+1)}\lla I^s/I^{s+1},
\end{equation}
for which $\pi_*X^{(s+1)}=R_*/I_*^{s+1}$ as an $R_*$-module.

If $M$ is an $R$-module which is an $E$ module spectrum,
Theorem~\ref{thm:UCSS-RQ-Conditions} provides a Universal
Coefficient Spectral Sequence
\[
\mathrm{E}_2^{*,*}=
\Ext^{p,q}_{E_*}(E^R_*R/I^s,M_*)\Lra\mathcal{D}_R(R/I^s,M)^{p+q}.
\]
Since $E^R_*R/I^s$ is $E_*$-free, this spectral sequence collapses
to give
\[
\mathcal{D}_R(R/I^s,M)^*=\Hom^*_{E_*}(E^R_*R/I^s,M_*).
\]
In particular, for $M=I^s/I^{s+1}$,
\[
\mathcal{D}_R(R/I^s,I^s/I^{s+1})^n=\Hom^n_{E_*}(E^R_*R/I^s,I_*^s/I_*^{s+1}).
\]
By~\eqref{eqn:R/Is-R/I(s+1)-Ext1} and Theorem~\ref{thm:R/Is-Resn}, there is
an element
\[
\tilde\del^{(s)}_*\in
%\Hom_{\mathcal{D}_{R_*}}^0(E_*\oTimes{R_*}\mathbf{K}^{(s-1)}_{*,*},I_*^s/I_*^{s+1}[-1])
%=
\Hom_{E_*}^0(E^R_*R/I^s,I_*^s/I_*^{s+1}[-1])
=
\Hom_{E_*}^1(E^R_*R/I^s,I_*^s/I_*^{s+1}),
\]
corresponding to an element $\delta_s\:R/I^s\lra\Sigma I^s/I^{s+1}$ inducing
a fibre sequence as in \eqref{eqn:FibSeq(s+1)}. It still remains to verify that
$\pi_*X^{(s+1)}=R_*/I_*^{s+1}$ as an $R_*$-module.

For this, we will use the resolutions $\mathbf{K}^{(s-1)}_{*,*}\lra R_*/I_*^s\ra0$
and $\mathbf{K}_{*,*}\lra R_*/I_*\ra0$. These free resolutions give rise to
cell $R$-module structures on $R/I^s$ and $E$. By~\cite{EKMM}, the $R$-module
$E\SSmash{R}R/I^s$ admits a cell structure with cells in one-one correspondence
with the elements of the obvious tensor product basis of
$\mathbf{K}_{*,*}\oTimes{R_*}\mathbf{K}^{(s-1)}_{*,*}$. Hence there is a resolution
by free $R_*$-modules
\[
\mathbf{K}_{*,*}\oTimes{R_*}\mathbf{K}^{(s-1)}_{*,*}\lra E^R_*R/I^s\ra0.
\]
There are morphisms of chain complexes
\[
\mathbf{K}^{(s-1)}_{*,*}
\xrightarrow{\rho_s}
\mathbf{K}_{*,*}\oTimes{R_*}\mathbf{K}^{(s-1)}_{*,*}
\xrightarrow{\tilde\delta_s}
\mathbf{Q}^{(s)}_{*,*}[-1],
\]
where $\rho_s$ is the obvious inclusion and $\tilde\delta_s$ is
a chain map lifting $\tilde\del^{(s)}_*$ which can be chosen so
that
\[
\tilde\delta_s(e_i\otimes x)=0.
\]
The effect of the composite $\tilde\delta_s\rho_s$ on the generator
$e_i\tilde u_{(j_1,\ldots,j_{s-1})}\in\mathbf{K}^{(s-1)}_{1,*}$ turns
out to be
\[
\tilde\del^{(s)}_*e_i\tilde u_{(j_1,\ldots,j_{s-1})}=
\tilde u_{(i,j_1,\ldots,j_{s-1})},
\]
while the elements of form $e_i\otimes\tilde u_{(j_1,\ldots,j_{k-1})}$ with
$k<s$ are annihilated. The composite homomorphism
\[
\mathbf{K}^{(s-1)}_{1,*}\xrightarrow{\tilde\delta_s\rho_s}
\mathbf{Q}^{(s)}_{0,*}[-1]\xrightarrow{\epsilon_1} I_*^s/I_*^{s+1}[-1]
\]
is a cocycle. There is a morphism of exact sequences
\[
\begin{CD}
0@<<<R_*/I_*^s @<<<\mathbf{K}^{(s-1)}_{0,*}@<<<\mathbf{K}^{(s-1)}_{1,*}
@<<<\mathbf{K}^{(s-1)}_{2,*}\\
@.    @|   @V\alpha_0 VV @V\alpha_1 VV @VVV  \\
0@<<<R_*/I_*^s  @<<< R_*/I_*^{s+1} @<<<I_*^s/I_*^{s+1} @<<< 0 @.
\end{CD}
\]
where the cohomology class
\[
[\alpha_1]\in\Ext_{R_*}^{1,*}(R_*/I_*^s,I_*^s/I_*^{s+1})
\]
represents the extension of $R_*$-modules on the bottom row. It is
easy to see that $[\alpha_1]=[\epsilon_1\tilde\delta_s\rho_s]$,
hence this class also represents the extension of $R_*$-modules
\[
0\la R_*/I_*^s \lla \pi_*X^{s+1}\lla I_*^s/I_*^{s+1}\la 0.
\]

There is a diagram of cofibre triangles
\[
\begin{xy}
\xymatrix{
R/I^{s+1}\ar[d]&I/I^{s+1}\ar[d]\ar[l]
%&I^2/I^{s+1}\ar[d]\ar[l]
&\ar[l]
&\cdots
&&I^{s-1}/I^{s+1}\ar[l]\ar[d]&I^s/I^{s+1}\ar[l]\ar[d]^{=}\\
R/I\ar[ur]&I/I^2\ar[ur]
%&I^2/I^3\ar[ur]
&&  &\ar[ur]&I^{s-1}/I^{s}\ar[ur]&I^s/I^{s+1}
}
\end{xy}
\]
and applying $E^R_*(\ )$ we obtain a spectral sequence converging to
$E^R_*R/I^{s+1}$ whose $\E_2$-term is the homology of the complex
\[
0\ra E^R_*R/I\xrightarrow{\del^{(1)}_*}E^R_*I/I^2\xrightarrow{\del^{(2)}_*}
E^R_*I^2/I^3\lra\cdots\xrightarrow{\del^{(s)}_*}E^R_*I^s/I^{s+1}\ra0,
\]
where the $\del^{(k)}_*$ are essentially the maps used to compute
$\Tor^{R_*}_*(R_*/I_*,R_*/I_*^{s+1})$ in~\cite{AB:HomRegQuot}. By
Proposition~\ref{prop:QComp-Exactness} and Theorem~\ref{thm:Tor(R/I,R/Is)},
this complex is exact except at the ends, where we have $\ker\del^{(1)}_*=E_*$.
As a result, this spectral sequence collapses at $\E_3$ giving the desired form
for $E^R_*R/I^{s+1}$.
\end{proof}
\begin{cor}\label{cor:R/Is-E*Proj}
For any $E$-module spectrum $N$ and $s\geq1$,
\[
N_R^*R/I^s\iso\Hom_{E_*}(E^R_*R/I^s,N_*).
\]
\end{cor}
\begin{proof}
This follows from Theorem~\ref{thm:R/I^s}(i).
\end{proof}
We will also use the following result.
\begin{cor}\label{cor:ER*R/Is-TowerSurj}
For $s\geq1$, the natural map
\begin{align*}
E^R_*R/I^{s+1}\lra E^R_*R/I^s,
\end{align*}
has image equal to $E_*=E^R_*R$.
\end{cor}
\begin{proof}
This follows from Theorem~\ref{thm:R/I^s}(ii).
\end{proof}
\begin{cor}\label{cor:R/Is-TowerSurj}
For any $E$-module spectrum $N$ and $s\geq1$,
\[
\colim_sN_R^*R/I^s\iso N_R^*R\iso N_*.
\]
\end{cor}
\begin{proof}
This is immediate from Corollaries~\ref{cor:R/Is-E*Proj}
and \ref{cor:ER*R/Is-TowerSurj} since
\[
\colim_s\Hom_{E_*}(E^R_*R/I^s,N_*)\iso\Hom_{E_*}(E_*,N_*).
\]
\end{proof}

\section{The $I$-adic tower and Adams Spectral Sequence}
\label{sec:I-adictower-ASS}

Continuing with the notation of Section~\ref{sec:I-adictower},
the first substantial result of this section is
\begin{thm}\label{thm:I-adic-holim}
The $I$-adic tower
\[
R/I\lla R/I^2\lla\cdots\lla R/I^s\lla R/I^{s+1}\lla\cdots
\]
has homotopy limit
\[
\holim_sR/I^s\simeq\hat{\mathrm{L}}^R_ES_R.
\]
\end{thm}
Our approach follows ideas of Bousfield~\cite{Bousfield:ASS}
where it is shown that the following Lemma implies
Theorem~\ref{thm:I-adic-holim}.
\begin{lem}\label{lem:E-nilpotentTow-Conditions}
Let $E=R/I$. Then the following are true. \\
{\rm i)} Each $R/I^s$ is $E$-nilpotent. \\
{\rm ii)} For each $E$-nilpotent $R$-module $M$,
\[
\colim_s\mathcal{D}_R(R/I^s,M)^*=M_{-*}.
\]
\end{lem}
\begin{proof}
(i) is proved by an easy induction on $s\geq1$. \\
(ii) is a consequence of Corollary~\ref{cor:R/Is-TowerSurj}.
\end{proof}

Since the maps $R_*/I_*^{s+1}\lra R_*/I_*^s$ are surjective,
from the standard exact sequence for $\pi_*(\ )$ of a homotopy
limit we have
\begin{equation}\label{eqn:LQMU-nilcomp-pi*}
\pi_*\hat{\mathrm{L}}^R_ES_R=\lim_sR_*/I_*^s.
\end{equation}
We can generalize this to the case where $E$ is a commutative
localized regular quotient.
\begin{thm}\label{thm:LQMU-nilcomp}
Let $E=R/I[X^{-1}]$ be a commutative localized regular quotient
of $R$. Then
\[
\pi_*\hat{\mathrm{L}}^R_ES_R=R_*[X^{-1}]\sphat_{I_*}
=\invlim_sR_*[X^{-1}]/I_*^s.
\]
%%%replaced
%The natural map $S_R\lra\hat{\mathrm{L}}^R_ES_R$ is an $E$-equivalence,
%%%by
If the regular sequence generating $I_*$ is finite, then the
natural map $S_R\lra\hat{\mathrm{L}}^R_ES_R$ is an
$E$-equivalence,
%%%
hence
\begin{align*}
\mathrm{L}^R_ES_R&\simeq\hat{\mathrm{L}}^R_ES_R, \\
\pi_*\mathrm{L}^R_ES_R&=R_*[X^{-1}]\sphat_{I_*}.
\end{align*}
\end{thm}
\begin{proof}
The first statement is easy to verify.

By Remark~\ref{rem:Ehat-Invce}, to simplify notation we may as well
replace $R$ by $R[X^{-1}]$ and so assume that $E=R/I$ is a commutative
regular quotient of $R$.

Using the Koszul complex $(\Lambda_{R_*}(e_j:j),\d)$, we see that
$\Tor^{R_*}_{*,*}(E_*,(R_*)\sphat_{I_*})$ is the homology of the
complex
\[
\Lambda_{R_*}(e_j:j)\oTimes{R_*}(R_*)\sphat_{I_*}
=\Lambda_{(R_*)\sphat_{I_*}}(e_j:j)
\]
with differential $\d'=\d\otimes1$.
%%%replaced
%Since the sequence $u_j$ remains
%regular in $(R_*)\sphat_{I_*}$, this complex provides a free resolution
%of $E_*=R_*/I_*$ as an $(R_*)\sphat_{I_*}$-module.
%%%by
Since the sequence $u_j$ remains regular in $(R_*)\sphat_{I_*}$,
this complex provides a free resolution of $E_*=R_*/I_*$ as an
$(R_*)\sphat_{I_*}$-module (this is \emph{false} if the sequence
$u_j$ is not finite).
%%%
Hence we have
\[
\Tor^{R_*}_{*,*}(E_*,(R_*)\sphat_{I_*})
=
\Tor^{(R_*)\sphat_{I_*}}_{*,*}(E_*,(R_*)\sphat_{I_*})
=E_*.
\]

To calculate $E_*^R\hat{\mathrm{L}}^R_ES_R$ we may use the K\"unneth
Spectral Sequence of~\cite{EKMM},
\[
\E_2^{s,t}=\Tor^{R_*}_{s,t}(E_*,\hat{\mathrm{L}}^R_ES_R)
\Lra E^R_{s+t}\hat{\mathrm{L}}^R_ES_R.
\]
%%%replaced
%By Theorem~\ref{thm:LQMU-nilcomp},
%%%by
By the first part,
%%%
the $\E_2$-term is
\[
\Tor^{R_*}_{*,*}(E_*,(R_*)\sphat_{I_*})=E_*=E^R_*R.
\]
Hence the natural homomorphism
\[
E^R_*S_R\lra E^R_*\hat{\mathrm{L}}^R_ES_R
\]
is an isomorphism.
\end{proof}
%%%replaced
%A generalization of this result
%\begin{thm}\label{thm:LQ-R}
%Let $E$ be a commutative localized regular quotient of $R$.
%Then for any finite cell $R$-module $M$, the natural map
%$M\lra\hat{\mathrm{L}}^R_EM$ is an $E$-equivalence, hence
%\begin{align*}
%\mathrm{L}^R_EM&\simeq\hat{\mathrm{L}}^R_EM, \\
%\pi_*\mathrm{L}^R_EM&=M_*[X^{-1}]\sphat_{I_*}
%=R_*[X^{-1}]\sphat_{I_*}\oTimes{R_*}M_*.
%\end{align*}
%\end{thm}
%%%by

If the sequence $u_j$ is infinite, the calculation of this proof
shows that
\[
E^R_*\hat{\mathrm{L}}^R_ES_R=(R_*)\sphat_{I_*}/I_*\neq
R_*/I_*=E_*S_R
\]
and the Adams Spectral Sequence does not converge to the homotopy
of the $E$-localization.

An induction on the number of cells of $M$ proves a generalization
of Theorem~\ref{thm:LQMU-nilcomp}.
\begin{thm}\label{thm:LQ-R}
Let $E$ be a commutative localized regular quotient of $R$ and $M$
a finite cell $R$-module. Then
\[
\pi_*\hat{\mathrm{L}}^R_EM=M_*[X^{-1}]\sphat_{I_*}
=R_*[X^{-1}]\sphat_{I_*}\oTimes{R_*}M_*.
\]
If the regular sequence generating $I_*$ is finite, then the
natural map $M\lra\hat{\mathrm{L}}^R_EM$ is an $E$-equivalence,
hence
\begin{align*}
\L^R_EM&\simeq\hat{\mathrm{L}}^R_EM, \\
\pi_*\L^R_EM&=M_*[X^{-1}]\sphat_{I_*}
=R_*[X^{-1}]\sphat_{I_*}\oTimes{R_*}M_*.
\end{align*}
\end{thm}
%%%

The reader may wonder if the following conjecture is true, the algebraic
issue being that it does not appear to be true that for a commutative ring
$A$, the extension $A\lra A\sphat_J$ is always flat for an ideal $J\ideal A$,
a Noetherian condition normally being required to establish such a result.
\begin{conj}\label{conj:LQ}
The conclusion of Theorem~\ref{thm:LQ-R} holds when $E$ is any commutative
localized quotient of $R$.
\end{conj}

\section{Some examples associated with $\MU$}
\label{sec:LRQ-MU}

An obvious source of commutative localized regular quotients is the
commutative $S$-algebra $R=\MU$ and we will describe some important
examples. It would appear to be algebraically simpler to work with
$\BP$ at a prime $p$ in place of $\MU$, but at the time of writing,
it seems not to be known whether $\BP$ admits a commutative $S$-algebra
structure.

\subsection*{Example A: $\MU\lra H\F_p$.}

Let $p$ be a prime. By considering the Eil\-en\-berg\--Mac~Lane
spectrum $H\F_p$ as a commutative $\MU$-algebra~\cite{EKMM},
we can form $H\F_p\SSmash{\MU}H\F_p$. The K\"unneth Spectral
Sequence gives
\[
\E^2_{s,t}=\Tor_{s,t}^{\MU_*}(\F_p,\F_p)
\Lra
{H\F_p\,}^{\MU}_{s+t}H\F_p.
\]
Using a Koszul complex over $\MU_*$, it is straightforward to see that
\[
\E^2_{*,*}=\Lambda_{\F_p}(\tau_j:j\geq0),
\]
the exterior algebra over $\F_p$ with generators
$\tau_j\in\E^2_{1,2j}$.

Taking $R=\MU$ and $E=H\F_p$, we obtain a spectral sequence
\[
\E_2^{s,t}(\MU)=
\Coext_{\Lambda_{\F_p}(\tau_j:j\geq0)}^{s,t}(\F_p,\F_p)
\Lra
\pi_{s+t}
%%%replaced
%\L
%%%by
\hat{\L}
%%%
^\MU_{H\F_p}S_\MU,
\]
where $I_\infty\ideal\MU_*$ is generated by $p$ together with
all positive degree elements, so $\MU_*/I_\infty=\F_p$. Also,
\[
\pi_*
%%%replaced
%\L
%%%by
\hat{\L}
%%%
^\MU_{H\F_p}S_\MU=(\MU_*)\sphat_{I_\infty}.
\]
More generally, for a finite cell $\MU$-module $M$, the Adams
Spectral Sequence has the form
\[
\E_2^{s,t}(M)=
\Coext_{\Lambda_{\F_p}(\tau_j:j\geq0)}^{s,t}(\F_p,{H\F_p}^\MU_*M)
\Lra
\pi_{s+t}
%%%replaced
%\L
%%%by
\hat{\L}
%%%
^\MU_{H\F_p}M,
\]
where
\[
\pi_*
%%%replaced
%\L
%%%by
\hat{\L}
%%%
^\MU_{H\F_p}M=(M_*)\sphat_{I_\infty}.
\]

\subsection*{Example B: $\MU\lra\En$.}

By~\cite{EKMM,Strickland:MU}, the Johnson-Wilson spectrum $\En$ at an
\emph{odd} prime $p$ is a commutative $\MU$-ring spectrum. According
to proposition~2.10 of~\cite{Strickland:MU}, at the prime $2$ a
certain modification of the usual construction also yields a commutative
$\MU$-ring spectrum which we will still denote by $\En$ rather than
Strickland's $\En'$. In all cases we can form the commutative $\MU$-ring
spectrum $\En\SSmash{\MU}\En$ and there is a K\"unneth Spectral Sequence
\[
\E^2_{s,t}=\Tor_{s,t}^{\MU_*}(\En_*,\En_*)
\Lra
\En^{\MU}_{s+t}\En.
\]
By using a Koszul complex for $\MUn_*$ over $\MU_*$ and localizing
at $v_n$, we find that
\[
\E^2_{*,*}=
\Lambda_{\En_*}(\tau_j:\text{$j\geq1$ and $j\neq p^k-1$ with $1\leq k\leq n$}),
\]
where $\Lambda$ denotes an exterior algebra and $\tau_j\in\E^2_{1,2j}$.
So
\[
\En^{\MU}_*\En=
\Lambda_{\En_*}(\tau_j:
\text{$j\geq1$ and $j\neq p^k-1$ with $1\leq k\leq n$})
\]
as an $\En_*$-algebra.

When $R=\MU$ and $E=\En$, we obtain a spectral sequence
\[
\E_2^{s,t}(\MU)=
\Coext_{\Lambda_{\En_*}(\tau_j:j\geq n+1)}^{s,t}(\En_*,\En_*)
\Lra
\pi_{s+t}
%%%replaced
%\L
%%%by
\hat{\L}
%%%
^\MU_{\En}\MU,
\]
where
\[
\pi_*
%%%replaced
%\L
%%%by
\hat{\L}
%%%
^\MU_{\En}\MU=(\MU_*)_{(p)}[v_n^{-1}]\sphat_{J_{n+1}}
\]
and
\[
J_{n+1}=(\ker\:(\MU_*)_{(p)}[v_n^{-1}]\lra\En_*)\ideal\MU_*[v_n^{-1}].
\]
In the $\E_2$-term we have
\[
\E_2^{s,t}(\MU)=
\En_*[U_j:\text{$0\leq j\neq p^k-1$ for $0\leq k\leq n$}],
\]
with generator $U_j\in\E_2^{1,2j+1}(\MU)$ corresponding to
an exterior generator in $\En^{\MU}_*\En$ associated with
a polynomial generator of $\MU_*$ in degree $2j$ lying in
$\ker\MU_*\lra\En_*$.

More generally, for a finite cell $\MU$-module $M$,
\[
\E_2^{s,t}(M)=
\Coext_{\Lambda_{\En_*}(\tau_j:j\geq n+1)}^{s,t}(\En_*,\En^{\MU}_{\ *}M),
\Lra
\pi_{s+t}
%%%replaced
%\L
%%%by
\hat{\L}
%%%
^\MU_{\En}M,
\]
where
\[
\pi_*
%%%replaced
%\L
%%%by
\hat{\L}
%%%
^\MU_{\En}M=M\sphat_{J_{n+1}}
=(\MU_*)_{(p)}[v_n^{-1}]\sphat_{J_{n+1}}\oTimes{\MU}M.
\]

\subsection*{Example C: $\MU\lra\Kn$.}

We know from~\cite{EKMM,Strickland:MU} that for an odd prime $p$, the
spectrum $\Kn$ representing the $n$\,th Morava $K$-theory $\Kn^*(\ )$
is a commutative $\MU$ ring spectrum. There is a K\"unneth Spectral
Sequence
\[
\E^2_{s,t}=\Tor_{s,t}^{\MU_*}(\Kn_*,\Kn_*)
\Lra
\Kn^{\MU}_{s+t}\Kn,
\]
and we have
\[
\E^2_{*,*}=\Lambda_{\Kn_*}(\tau_j:0\leq j\neq p^n-1).
\]

Taking $R=\MU$ and $E=\Kn$, we obtain a spectral sequence
\[
\E_2^{s,t}(\MU)=
\Coext_{\Lambda_{\Kn_*}(\tau_j:\text{$0\leq j\neq n$})}^{s,t}(\Kn_*,\Kn_*)
\Lra\pi_{s+t}
%%%replaced
%\L
%%%by
\hat{\L}
%%%
^{\MU}_{\Kn}\MU,
\]
where
\[
\pi_*
%%%replaced
%\L
%%%by
\hat{\L}
%%%
^{\MU}_{\Kn}\MU=(\MU_*)\sphat_{I_{n,\infty}}
\]
with $I_{n,\infty}=\ker\MU_*\lra\Kn_*$. In the $\E_2$-term
we have
\[
\E_2^{s,t}(\MU)=
\En_*[U_j:\text{$0\leq j\neq p^n-1$}],
\]
with generator $U_j\in\E_2^{1,2j+1}(\MU)$ corresponding to
an exterior generator in $\En^{\MU}_*\En$ associated with
a polynomial generator of $\MU_*$ in degree $2k$ lying in
$\ker\MU_*\lra\En_*$ (or when $j=0$, associated with $p$).

More generally, for a finite cell $\MU$-module $M$,
\[
\E_2^{s,t}(M)=
\Coext_{\Lambda_{\Kn_*}(\tau_j:\text{$0\leq j\neq n$})}^{s,t}(\Kn_*,\Kn^{\MU}_*M)
\Lra\pi_{s+t}
%%%replaced
%\L
%%%by
\hat{\L}
%%%
^{\MU}_{\Kn}M,
\]
where
\[
\pi_*
%%%replaced
%\L
%%%by
\hat{\L}
%%%
^{\MU}_{\Kn}M=(M_*)_{I_n,\infty}\sphat
=(\MU_*)\sphat_{I_{n,\infty}}\oTimes{\MU_*}M_*.
\]

\section*{Concluding remarks}

There are several outstanding issues raised by our work.

Apart from the question of whether it is possible to weaken the assumptions
from (commutative) regular quotients to a more general class, it seems reasonable
to ask whether the internal $I$-adic tower is one of $R$ ring spectra.
%%%replaced
%Since $\L^R_ER=\ds\holim_s R/I^s$,
%%%by
Since $\L^R_ER=\ds\holim_s R/I^s$ (at least when $I_*$ is finitely
generated),
%%%
the localization theory of~\cite{EKMM,Wolbert} shows that this can be
realized as a commutative $R$-algebra. However, showing that each $R/I^s$
is an $R$ ring spectrum or even an $R$-algebra seem to involve far more
intricate calculations. We expect that this will turn out to be true and even
that the tower is one of $R$-algebras. This should involve techniques similar
to those of~\cite{Lazarev,AB+AJ:Brave-MU}. It is also worth noting that our
proofs make no distinction between the cases where $I_*\ideal R_*$ is infinitely
or finitely generated. There are a number of algebraic simplifications possible
in the latter case, however we have avoided using them since the most interesting
examples we know are associated with infinitely generated regular ideals in
$\MU_*$. The spectra $E_n$ of Hopkins, Miller \emph{et al.} have Noetherian
homotopy rings and there are towers based on powers of their maximal ideals
similar to those in the first author's previous work~\cite{AB:Ainfty,AB-UW:Bockop}.

We also hope that our preliminary exploration of Adams Spectral Sequences for
$R$-modules will lead to further work on this topic, particularly in the case
$R=\MU$ and related examples. A more ambitious project would be to investigate
the commutative $S$-algebra $\MSp$ from this point of view, perhaps reworking
the results of Vershinin, Gorbounov and Botvinnik in the context of
$\MSp$-modules~\cite{Botvinnik,Vershinin}.

\Addresses\recd

\end{document}